\documentclass[reqno]{amsart}

\usepackage{amssymb}
\theoremstyle{plain}
\newtheorem{theorem}{Theorem}[section]
\newtheorem{lemma}[theorem]{Lemma}

\newtheorem{proposition}[theorem]{Proposition}
\newtheorem{problem}[theorem]{Problem}
\theoremstyle{remark}

\newcommand\QQ{\mathcal{Q}}
\newcommand\sbl[1]{\langle#1\rangle}   % subloop or subquasigroup
\newcommand\ld{\backslash} % left division
\newcommand\rd{/}  % right division
\newcommand\opp{\mathrm{op}}

\title[Axioms for Trimedial Quasigroups]
{Axioms for Trimedial Quasigroups}

\author[M.~K.~Kinyon]{Michael~K.~Kinyon}
\address{Department of Mathematical Sciences \\
Indiana University South Bend \\
South Bend, IN 46634 USA}
\email{mkinyon@iusb.edu}
\urladdr{http://mypages.iusb.edu/\symbol{126}mkinyon}
\author[J.~D.~Phillips]{J.~D.~Phillips}
\address{Department of Mathematics \& Computer Science \\
Wabash College \\
Crawfordsville, IN 47933 USA}
\email{phillipj@wabash.edu}
\urladdr{http://www.wabash.edu/depart/math/faculty.html{\#}Phillips}

\date{\today}

\subjclass{20N05}
\keywords{trimedial quasigroup, F-quasigroup, semimedial quasigroup}

\newcommand{\lemref}[1]{Lemma~\ref{#1}}
\newcommand{\thmref}[1]{Theorem~\ref{#1}}
\newcommand{\propref}[1]{Proposition~\ref{#1}}

\begin{document}

\begin{abstract}
We give new equations that axiomatize the variety of trimedial quasigroups.
We also improve a standard characterization by showing that right
semimedial, left F-quasigroups are trimedial.
\end{abstract}

\maketitle

%%%%%%%%%%
\section{Introduction}
\label{sec:intro}

A quasigroup $\QQ = (Q;\;\cdot, \ld, \rd)$ is a set $Q$ with three binary operations
$\cdot, \ld, \rd : Q\times Q\to Q$ satisfying the equations:
\[
\begin{array}{c}
x\ld (x\cdot y) = y \qquad\qquad (x\cdot y) \rd y = x\\
x\cdot (x\ld y) = y  \qquad\qquad (x\rd y)\cdot y = x
\end{array}
\]
Basic references for quasigroup theory are \cite{Bel}, \cite{Br}, \cite{CPS}, \cite{Pf}.

A quasigroup is \emph{medial} if it satisfies the identity
\begin{align*}
xy\cdot uv &= xu\cdot yv \tag{$M$}
\end{align*}
A quasigroup is \emph{trimedial} if every subquasigroup generated by three
elements is medial. Medial quasigroups have also been called abelian, entropic,
and other names, while trimedial quasigroups have also been called triabelian,
terentropic, etc. (See Chap. IV of \cite{CPS}, especially p. 120, for further background.)
The classic Toyoda-Bruck theorem asserts that every medial quasigroup is isotopic
to an abelian group \cite{To} \cite{Br44}. This result was generalized by Kepka to
trimedial quasigroups: every trimedial quasigroup is isotopic to a commutative
Moufang loop \cite{Ke76}.

There are two distinct, but related, generalizations of trimedial quasigroups. The
variety of \emph{semimedial} quasigroups (also known as weakly abelian, weakly
medial, etc.) is defined by the equations
\begin{align*}
  xx\cdot yz &= xy\cdot xz \tag{$S_l$} \\
  zy\cdot xx &= zx\cdot yx  \tag{$S_r$}
\end{align*}
A quasigroup satisfying ($S_l$) (resp. ($S_r$)) is said to be
\emph{left} (resp. \emph{right}) \emph{semimedial}. 
Every semimedial quasigroup is isotopic to a commutative Moufang loop \cite{Ke76}.
(In the trimedial case, the isotopy has a more restrictive form; see the cited references for
details.)

The variety of \emph{F-quasigroups} was introduced by Murdoch in \cite{Mu},
the same paper in which he introduced what we now call medial quasigroups. F-quasigroups
are defined by the equations
\begin{align*}
  x \cdot yz &= xy\cdot (x\ld x)z  \tag{$F_l$} \\
  zy\cdot x &= z(x\rd x) \cdot yx \tag{$F_r$} 
\end{align*}
A quasigroup satisfying ($F_l$) (resp. ($F_r$)) is said to be a
\emph{left} (resp. \emph{right}) \emph{F-quasigroup}. 
Murdoch did not actually name the variety of F-quasigroups. We thank one
of the referees for suggesting that the earliest use of the name might be in \cite{BO}.

One among many links between these two generalizations of trimedial quasigroups
is the following (\cite{Ke79}, Prop. 6.2).

\begin{proposition}
\label{prop:tri=F+semi}
A quasigroup is trimedial if and only if it is a semimedial, left (or right) F-quasigroup.
\end{proposition}

Together with Kepka, we have been investigating the structure of F-quasigroups, and
have shown that every loop isotopic to an F-quasigroup is Moufang. Full details will
appear elsewhere \cite{KKP}. As part of that investigation, we were led to consider
the following equations, which are similar in form to ($F_l$), ($F_r$):
\begin{align*}
x \cdot yz &= (x\rd x) y \cdot xz \tag{$E_l$} \\
zy\cdot x &= zx\cdot y (x\ld x) \tag{$E_r$}
\end{align*}

The main result of the present paper is the following.

\begin{theorem}
\label{thm:main}
A quasigroup is trimedial if and only if it satisfies ($E_l$) and ($E_r$).
\end{theorem}

Kepka \cite{Ke76} \cite{Ke78} showed that the variety of trimedial quasigroups is
axiomatized by the semimedial laws ($S_l$), ($S_r$), and by the equation
$(x\cdot xx)\cdot uv = xu \cdot (xx\cdot v)$. Later \cite{KP} we showed
that ($S_l$) is redundant. \thmref{thm:main} offers a more symmetric
alternative.

As an auxiliary result, we will also use ($E_l$) and ($E_r$)
to obtain the following improvement of \propref{prop:tri=F+semi}.

\begin{theorem}
\label{thm:main2}
Let $\QQ$ be a quasigroup. The following are equivalent:
\begin{enumerate}
\item $\QQ$ is trimedial.
\item $\QQ$ is a right semimedial, left F-quasigroup.
\item $\QQ$ is a left semimedial, right F-quasigroup.
\end{enumerate}
\end{theorem}

Our investigations were aided by the equational reasoning tool OTTER
developed by McCune \cite{Mc}. We thank T.~Kepka for suggesting
that ($E_l$), ($E_r$) might axiomatize an interesting variety
of quasigroups; he was certainly correct.

\section{Proofs}

Our strategy for proving \thmref{thm:main} is to use \propref{prop:tri=F+semi}: we
will show that a quasigroup satisfying ($E_l$), ($E_r$) is a semimedial,
F-quasigroup. First we introduce some notation for local right and left unit elements
in a quasigroup:
\[
e(x) := x \ld x  \qquad  f(x) := x\rd x
\]

If $\QQ = (Q;\cdot, \ld, \rd)$ is a quasigroup, then so are the
\emph{left parastrophe} $\QQ_l := (Q; \ld, \cdot, \rd_{\opp})$,
the \emph{right parastrophe} $\QQ_r := (Q; \rd, \cdot, \ld_{\opp})$,
and the \emph{opposite parastrophe}
$\QQ_{\opp} := (Q; \cdot_{\opp}, \rd, \ld)$, 
where for a binary operation $*$, we use $*_{\opp}$ to denote the
opposite operation. Note that $(\QQ_l)_l = \QQ$, $(\QQ_r)_r = \QQ$,
$(\QQ_{\opp})_{\opp}= \QQ$, $(\QQ_{\opp})_l = (\QQ_r)_{\opp}$,
and $(\QQ_{\opp})_r = (\QQ_l)_{\opp}$.  For a more complete discussion
of parastrophes, including alternative notation conventions, see \cite{Bel}, \cite{Pf}.

We state the following obvious result formally for later ease of reference.

\begin{lemma}
\label{lem:obvious}
Let $\QQ = (Q;\cdot, \ld, \rd)$ be a quasigroup.
\begin{enumerate}
\item $\QQ$ satisfies ($F_l$) if and only if $\QQ_{\opp}$ satisfies ($F_r$).
\item $\QQ$ satisfies ($S_l$) if and only if $\QQ_{\opp}$ satisfies ($S_r$).
\item $\QQ$ satisfies ($E_l$) if and only if $\QQ_{\opp}$ satisfies ($E_r$).
\end{enumerate}
\end{lemma}

Parts (1) and (2) of the following lemma are well-known, although
the authors have not been able to find a specific reference.

\begin{lemma}
\label{lem:para}
Let $\QQ = (Q;\cdot, \ld, \rd)$ be a quasigroup.
\begin{enumerate}
\item $\QQ$ is a left F-quasigroup if and only if $\QQ_l$ is left semimedial.
\item $\QQ$ is a right F-quasigroup if and only if $\QQ_r$ is right semimedial.
\item $\QQ$ satisfies ($E_l$) if and only if $\QQ_l$ satisfies ($E_l$).
\item $\QQ$ satisfies ($E_r$) if and only if $\QQ_r$ satisfies ($E_r$).
\end{enumerate}
\end{lemma}

\begin{proof}
 For (1): In $\QQ_l$, ($S_l$) is $e(x)\ld (y\ld z) = (x\ld y)\ld (x\ld z)$.
Multiply on the left by $e(x)$, replace $y$ with $xy$, and $z$ with $xz$ to get
$(xy)\ld (xz) = e(x)(y\ld z)$. Now replace $z$ with $yz$ and multiply on the
left by $xy$ to get $x\cdot yz = xy\cdot e(x)z$, which is ($F_l$) in $\QQ$.
Since $(\QQ_l)_l = \QQ$, the converse also holds.

For (2): By part (1), ($F_l$) holds in $\QQ_{\opp}$ iff ($S_l$) holds in
$(\QQ_{\opp})_l = (\QQ_l)_{\opp}$. The desired result now follows from
\lemref{lem:obvious}.

For (3): In $\QQ_l$, ($E_l$) is $((x\rd_{\opp} x)\ld y)\ld (x\ld z) = x\ld (y\ld z)$.
Multiply on the left by $f(x)\ld y$, and replace $z$ with $yz$ to get
$x\ld (yz) = (f(x)\ld y)(x\ld z)$. Multiply on the left by $x$, replace $y$ with $f(x)y$,
and $z$ with $xz$ to get $f(x)y\cdot xz = x\cdot yz$, which is ($E_l$) in
$\QQ$. Since $(\QQ_l)_l = \QQ$, the converse also holds.

For (4): By part (3), ($E_l$) holds in $\QQ_{\opp}$ iff ($E_l$) holds in
$(\QQ_{\opp})_l = (\QQ_r)_{\opp}$.  The desired result now follows from
\lemref{lem:obvious}.
\end{proof}

\begin{lemma}
\label{lem:end}
Let $\QQ$ be a quasigroup.
\begin{enumerate}
\item If $\QQ$ satisfies ($E_l$), then $f: Q \to Q$ is an endomorphism of $\QQ$.
\item If $\QQ$ satisfies ($E_r$), then $e: Q \to Q$ is an endomorphism of $\QQ$.
\item If $\QQ$ is a left F-quasigroup, then $e: Q\to Q$ is an endomorphism of $\QQ$.
\item If $\QQ$ is a right F-quasigroup, then $f : Q\to Q$ is an endomorphism of $\QQ$.
\end{enumerate}
\end{lemma}

\begin{proof}
In each case, it is enough to show that the multiplication is preserved.

For (1): $f(x)f(y)\cdot xy = f(x)x\cdot f(y)y = xy$, and so
$f(x)f(y) = (xy)\rd (xy) = f(xy)$.

For (2): Since $f : Q\to Q$ is an endomorphism of $\QQ$ iff $e : Q\to Q$
is an endomorphism of $\QQ_{\opp}$, this follows from part (1) and
\lemref{lem:obvious}(3).

For (3): if ($F_l$) holds, then $xy\cdot e(x)e(y) = x\cdot ye(y) = xy$,
and so $e(x)e(y) = (xy)\ld (xy) = e(xy)$ (\cite[p. 38, eq. (32)]{BO}, \cite[Lemma 4.2]{Ke79}, \cite{BB}).

For (4): This follows from part (3) and \lemref{lem:obvious}(1).
\end{proof}

\begin{lemma}
\label{lem:commute}
Let $\QQ$ be a quasigroup. If $e:Q\to Q$ or $f:Q\to Q$ is an endomorphism
of $\QQ$, then $f(e(x)) = e(f(x))$ for all $x\in Q$.
\end{lemma}

\begin{proof}
If $f$ is an endomorphism, then $f(e(x)) = f(x)\ld f(x) = e(f(x))$, and the
case where $e$ is an endomorphism is similar.
\end{proof}

\begin{lemma}
\label{lem:E_l-equiv} Let $\QQ$ be a quasigroup.
\begin{enumerate}
\item If $\QQ$ satisfies ($E_l$), then $\QQ$ is a left F-quasigroup if and only
if it is left semimedial.
\item If $\QQ$ satisfies ($E_r$), then $\QQ$ is a right F-quasigroup if and
only if it is right semimedial.
\end{enumerate}
\end{lemma}

\begin{proof}
For (1): Assume $\QQ$ satisfies ($E_l$). By \lemref{lem:para}(3), 
($E_l$) holds in $\QQ_l$. We will prove the implication
($S_l$) $\implies$ ($F_l$). Since this will also hold in 
$\QQ_l$, it will follow from \lemref{lem:para}(1) that the implication
($F_l$) $\implies$ ($S_l$) will hold in $\QQ$. Now
if $\QQ$ is left semimedial, then 
\[
\begin{array}{rcll}
xx\cdot yz &=& xy\cdot xz & \text{by\ }(S_l) \\
&=& f(xy)x \cdot (xy\cdot z) & \text{by\ }(E_l) \\
&=& (f(x)f(y)\cdot x)\cdot (xy\cdot z) & \text{by \lemref{lem:end}} \\
&=& (x\cdot f(y)e(x))\cdot (xy\cdot z) & \text{by\ }(E_l) \\
&=& xx\cdot (f(y)e(x)\cdot (x\ld (xy\cdot z))) & \text{by\ }(S_l)
\end{array}
\]
Cancelling and replacing $z$ with $e(x)z$, we have
\[
\begin{array}{rcll}
f(y)e(x)\cdot (x \ld (xy\cdot e(x)z)) &=& y\cdot e(x)z & \\
&=& f(y)e(x)\cdot yz & \text{by\ }(E_l)
\end{array}
\]
Cancelling, we obtain $yz = x\ld (xy\cdot e(x)z)$ or 
$x\cdot yz = xy\cdot e(x)z$, which is ($F_l$).

For (2): If $\QQ$ satisfies ($E_r$), then $\QQ_{\opp}$ satisfies ($E_l$)
by \lemref{lem:obvious}(3). By part (1), ($F_l$) and ($S_l$) are equivalent
in $\QQ_{\opp}$, and so ($F_r$) and ($S_r$) are equivalent in $\QQ$ by
\lemref{lem:obvious}(1,2).
\end{proof}

\begin{lemma}
\label{lem:E-implies-F}
A quasigroup satisfying ($E_l$), ($E_r$) is an F-quasigroup.
\end{lemma}

\begin{proof}
We will show ($E_l$), ($E_r$) $\implies$ ($F_l$). Since this
implication will also hold in the opposite parastrophe (by \lemref{lem:obvious}(3)),
($F_r$) will follow from \lemref{lem:obvious}(1). Thus we compute 
\[
\begin{array}{rcll}
x\cdot yz &=& f(x)y\cdot xz & \text{by\ }(E_l) \\
&=& [f(x)e(f(x))\cdot x(x\ld y)]\cdot xz & \\
&=& [x\cdot e(f(x))(x\ld y)]\cdot xz & \text{by\ }(E_l) \\
&=& [x\cdot xz]\cdot [e(f(x))(x\ld y)\cdot e(xz)] & \text{by\ }(E_r) \\
&=& [x\cdot xz]\cdot [f(e(x))(x\ld y)\cdot e(x)e(z)] & \text{by\ Lemmas~\ref{lem:end}~and~\ref{lem:commute}} \\
&=& [x\cdot xz]\cdot [e(x)\cdot (x\ld y)e(z)] & \text{by\ }(E_l)
\end{array}
\]
Now $xz \cdot (x\ld y)e(z) = yz$ by ($E_r$), and so
\[
x\cdot yz = [x\cdot xz] \cdot e(x)[(xz)\ld (yz)]
\]
Now replace $y$ with $(xz\cdot y)\rd z$ to obtain
\[
x(xz\cdot y) = [x\cdot xz] \cdot e(x)y .
\]
Finally replace $z$ with $x\ld z$ to get
\[
x\cdot zy = xz\cdot e(x)y,
\]
which is ($F_l$).
\end{proof} 

We can now prove our main result.

\begin{proof}[Proof of \thmref{thm:main}]
Suppose $\QQ$ is trimedial. Since each of ($E_l$) and ($E_r$) is
a special case in three variables of the medial identity ($M$), these
identities will hold in $\QQ$. Indeed, fix $a,b,c \in Q$. Then the
subquasigroup $\sbl{a,b,c}$ generated by $\{ a,b,c\}$
is medial. Taking $x = f(a)$, $y = a$, $u = b$, $v = c$
in ($M$), we obtain $a\cdot bc = f(a)b\cdot ac$, while taking
$x = c$, $y = b$, $u = a$, $v = e(a)$ in ($M$) gives
$cb\cdot a = ca\cdot (b\cdot e(a))$. Since $a,b,c\in Q$ were chosen
arbitrarily, we have ($E_l$) and ($E_r$).

Conversely, if $\QQ$ satisfies ($E_l$) and ($E_r$), then by 
\lemref{lem:E-implies-F}, $\QQ$ is an F-quasigroup, and by
\lemref{lem:E_l-equiv}, $\QQ$ is semimedial. 
\propref{prop:tri=F+semi} complete the proof.
\end{proof}

\begin{lemma}
\label{lem:FSE}
\begin{enumerate}
\item A right semimedial, left F-quasigroup satisfies ($E_r$).
\item A left semimedial, right F-quasigroup satisfies ($E_l$).
\end{enumerate}
\end{lemma}

\begin{proof}
For (1): Suppose $\QQ$ satisfies ($S_r$) and ($F_l$). Then
\[
\begin{array}{rcll}
(xz\cdot y)\cdot e(xz)z &=& xz\cdot yz & \text{by\ }(F_l) \\
&=& xy\cdot zz & \text{by\ }(S_r) \\
&=& ((xy)\rd e(xz))e(xz)\cdot zz & \\
&=& ((xy)\rd e(xz))z \cdot e(xz)z & \text{by\ }(S_r)
\end{array}
\]
Cancelling and using \lemref{lem:end}(3), we have
$xz\cdot y = ((xy)\rd e(x)e(z))z$. Now $x(y\rd e(z))\cdot e(x)e(z) = xy$ by
($F_l$), and so $xz\cdot y = x(y\rd e(z))\cdot z$. Replacing $y$ with $ye(z)$,
we obtain $xz\cdot ye(z) = xy\cdot z$, which is ($E_r$).

For (2): If ($S_l$) and ($F_r$) hold in $\QQ$, then ($S_r$)
and ($F_l$) hold in $\QQ_{\opp}$ (\lemref{lem:obvious}(1,2)), and so
($E_r$) holds in $\QQ_{\opp}$ by part (1). By \lemref{lem:obvious}(3),
($E_l$) holds in $\QQ$.
\end{proof}

We now turn to our auxiliary result.

\begin{proof}[Proof of \thmref{thm:main2}]
Let $\QQ$ be a right semimedial, left F-quasigroup. By \lemref{lem:FSE}(1),
($E_r$) holds. By \lemref{lem:E_l-equiv}(2), $\QQ$ is an F-quasigroup.
We will now show that ($S_l$) holds. First, using ($S_r$), ($F_l$), and ($S_r$) again, we have
\[
((xy)\rd z) e(x)\cdot z^2 = xy\cdot e(x)z = x\cdot yz = (x\rd z)y \cdot z^2.
\]
Cancelling and dividing on the right by $e(x)$, we obtain
\begin{align*}
(xy)\rd z  = ((x\rd z)y) \rd e(x) . \tag{*}
\end{align*}
Next we use  ($S_r$), ($E_r$), and ($S_r$) again
to compute
\[
((xy)\rd e(z))z\cdot e(z)^2 = xy\cdot z = xz\cdot ye(z) = ((xz)\rd e(z))y\cdot e(z)^2 .
\]
Cancelling, we obtain
\begin{align*}
((xy)\rd e(z))z = ((xz)\rd e(z))y . \tag{**}
\end{align*}
Finally, we verify ($S_l$) as follows:
\[
\begin{array}{rcll}
xy\cdot xz &=& ((xy)\rd z)x \cdot z^2 & \text{by\ }(S_r) \\
&=& [((x\rd z)y) \rd e(x)]x \cdot z^2 & \text{by\ }(*) \\
&=& [((x\rd z)x) \rd e(x)]y \cdot z^2 & \text{by\ }(**) \\
&=& (x^2 \rd z)y \cdot z^2 & \text{by\ }(*) \\
&=& x^2 \cdot yz & \text{by\ }(S_r)
\end{array}
\]
Since we have shown that $\QQ$ is a semimedial, F-quasigroup, it
follows from \propref{prop:tri=F+semi} that $\QQ$ is trimedial.

On the other hand, if $\QQ$ is a right semimedial, left F-quasigroup, then
$\QQ_{\opp}$ satisfies ($S_l$) and ($F_r$) by \lemref{lem:obvious}(1,2).
By the preceding argument, $\QQ_{\opp}$ is a semimedial, F-quasigroup,
and thus so is $\QQ$ by \lemref{lem:obvious}(1,2). Once again,
\propref{prop:tri=F+semi} completes the proof.
\end{proof}

In closing, we note that further investigations suggest themselves. For example,
it would be of interest to determine the structure of quasigroups that are
only assumed to satisfy ($E_l$), or, in view of \lemref{lem:E_l-equiv},
those satisfying ($E_l$), ($S_l$) and ($F_l$). In this line
we pose a couple of problems.

\begin{problem}
\begin{enumerate}
\item Characterize the loop isotopes of quasigroups satisfying ($E_l$).
\item Characterize the loop isotopes of quasigroups satisfying ($E_l$),
($S_l$), and ($F_l$).
\end{enumerate}
\end{problem}

%%%%%%%%


\begin{thebibliography}{99}

\bibitem{Bel} V.~D.~Belousov,
\textit{Foundations of the Theory of Quasigroups and Loops},
Izdat. Nauka, Moscow, 1967 (Russian).
MR 36{\#}1569, Zbl 163:01801.

\bibitem{BO} V.~D.~Belousov,
About one quasigroup class. (Russian)
\textit{Uchenye zapiski Beltskogo gospedinstituta im. A. Russo},
\textbf{5} (1960), 29--44.

\bibitem{BB} V.~D.~Belousov,
\textit{Elements of Quasigroup Theory: A Special Course}
(in Russian), Kishinev State University Press, Kishinev, 1981.

\bibitem{Br44} R.~H.~Bruck,
Some results in the theory of quasigroups.
\textit{Trans. Amer. Math. Soc.} \textbf{55} (1944), 19--52.
MR 5,229d, Zbl 0063.00635.

\bibitem{Br} R.~H.~Bruck,
\textit{A Survey of Binary Systems},
Springer-Verlag, 1971.
MR 20{\#}76, Zbl 206:30301.

\bibitem{CPS} O.~Chein, H.~O.~Pflugfelder, and J.~D.~H.~Smith (eds.),
\textit{Quasigroups and Loops: Theory and Applications},
Sigma Series in Pure Math. \textbf{9}, Heldermann Verlag, 1990.
MR 93g:21033, Zbl 0719.20036.

\bibitem{Ke76} T.~Kepka,
Structure of triabelian quasigroups,
\textit{Comment. Math. Univ. Carolinae} \textbf{17} (1976), 229--240.
MR 53{\#}10965,  Zbl 0338.20097.

\bibitem{Ke78} T.~Kepka,
A note on WA-quasigroups,
\textit{Acta Univ. Carolin. Math. Phys.} \textbf{19} (1978), no. 2, 61--62.
MR 80b:20094, Zbl 0382.20057.

\bibitem{Ke79} T.~Kepka,
F-quasigroups isotopic to Moufang loops,
\textit{Czechoslovak Math. J.} \textbf{29} (1979), 62--83.
MR 80b:20095, Zbl 0444.20067.

\bibitem{KKP} T.~Kepka, M.~K.~Kinyon, and J.~D.~Phillips,
The structure of F-quasigroups, in preparation.

\bibitem{KP} M.~K.~Kinyon and J.~D.~Phillips,
A note on trimedial quasigroups,
\textit{Quasigroups \& Related Systems} \textbf{9} (2002), 65--66.
CMP 1 943 753,  Zbl pre01894750.

\bibitem{Mc} W.~W.~McCune,
\textit{OTTER\ 3.3 Reference Manual and Guide},
Technical Memorandum ANL/MCS-TM-263, Argonne National Laboratory, 2003;
or see: {\tt{http://www.mcs.anl.gov/AR/otter/}}

\bibitem{Mu} D.~C.~Murdoch,
Quasi-groups which satisfy certain generalized associative laws,
\emph{Amer. J. Math.} \textbf{61} (1939), 509--522.
 Zbl 0020.34702.

\bibitem{Pf} H.~O.~Pflugfelder,
\textit{Quasigroups and Loops: Introduction},
Sigma Series in Pure Math. \textbf{8}, Heldermann Verlag, Berlin, 1990.
MR 93g:20132, Zbl 719:20036.

\bibitem{To} K.~Toyoda,
On axioms of linear functions,
\textit{Proc. Imp. Acad. Tokyo} \textbf{17} (1941), 221--227.
MR 7,241g, Zbl 0061.02403.

\end{thebibliography}
\end{document}